\documentclass[authoryear,preprint,review,12pt]{elsarticle}
\usepackage{amsmath, amssymb, enumerate, amsthm}
\usepackage{bm}

\usepackage{geometry}
\usepackage{xcolor}

\usepackage[colorlinks,allcolors=blue]{hyperref}

\usepackage{booktabs}%for \toprule etc

\newcommand{\bigO}{\mathcal{O}_p}

\newcommand{\R}{\mathbb{R}}
\newcommand{\N}{\mathbb{N}}

\newcommand{\E}{\mathbb{E}}

\newcommand{\Ab}{\bm{A}}

\newcommand{\card}{{\rm card}}
\newcommand{\tr}{{\rm tr}}

\newcommand{\ind}{\perp \!\!\! \perp}

\newcommand{\xb}{\bm{x}}
\newcommand{\xbtilde}{\tilde{\bm{x}}}
\newcommand{\Xb}{\bm{X}}
\newcommand{\Xbtilde}{\tilde{\bm{X}}}
\newcommand{\Xtilde}{\tilde{X}}

\newcommand{\betab}{\bm{\beta}}
\newcommand{\betabhat}{\hat{\bm{\beta}}}
\newcommand{\betabtilde}{\tilde{\bm{\beta}}}
\newcommand{\betabtildehat}{\hat{\tilde{\bm{\beta}}}}
\newcommand{\betatildehat}{\hat{\tilde{\beta}}}

\newcommand{\betahat}{\hat{\beta}}

\newcommand{\ub}{\bm{u}}

\newcommand{\ubtilde}{\tilde{\bm{u}}}

\newcommand{\vb}{\bm{v}}

\newcommand{\vbtilde}{\tilde{\bm{v}}}
\newcommand{\vbtildehat}{\hat{\tilde{\bm{v}}}}

\newcommand{\vtilde}{\tilde{v}}

\newcommand{\lambdatilde}{\tilde{\lambda}}
\newcommand{\lambdatildehat}{\hat{\tilde{\lambda}}}

\newcommand{\Sigmab}{\bm{\Sigma}}

\newcommand{\Sigmabtilde}{\tilde{\Sigmab}}

\newcommand{\Sigmabhat}{\hat{\Sigmab}}

\newcommand{\Sigmabtildehat}{\hat{\tilde{\Sigmab}}}

\newcommand{\sigmatilde}{\tilde{\sigma}}

\newcommand{\Etab}{\bm{\mathrm{H}}}
\newcommand{\Etabtilde}{\tilde{\Etab}{}}

\newcommand{\Epsilonb}{\bm{\mathrm{E}}}

\newcommand{\bK}{\mathcal{K}}

\newcommand{\bP}{\mathcal{P}}
\newcommand{\bD}{\mathcal{D}}

\newcommand{\bDelta}{\mathit{\Delta}}

\newcommand{\dstar}{{d^\ast}}
\newcommand{\deltastar}{{\delta^\ast}}

\numberwithin{equation}{section} %number equation according to section

\newcommand{\keywords}[1]{\emph{Keywords:} #1}

\newcommand{\mathsubjclassSingle}[2]{\emph{#1:} #2}

\newtheorem{lemma}{Lemma}[section]
\newtheorem{theorem}{Theorem}[section]
\newtheorem{corollary}{Corollary}[section]

\newtheorem*{assumption*}{\assumptionnumber}
\providecommand{\assumptionnumber}{}
\makeatletter
\newenvironment{assumption}[2]
 {%
  \renewcommand{\assumptionnumber}{Assumption #1$^{#2}$}%
  \begin{assumption*}%
  \protected@edef\@currentlabel{#1}%
 }
 {%
  \end{assumption*}
 }
\makeatother

\makeatletter
\def\ps@pprintTitle{%
 \let\@oddhead\@empty
 \let\@evenhead\@empty
 \def\@oddfoot{\textit{Preprint}}%
 \let\@evenfoot\@oddfoot}
\makeatother

\begin{document}
\begin{frontmatter}

\author[1,2]{Jan O. Bauer}
\ead{j.bauer@unibas.ch}

\author[2]{Bernhard Drabant}
\ead{bernhard.drabant@dhbw-mannheim.de}

\address[1]{Faculty of Business and Economics,
University of Basel \\ Peter Merian-Weg 6, 4002 Basel, Switzerland.}

\address[2]{Baden-Wuerttemberg Cooperative State University Mannheim,
Coblitzallee 1-9, 68163 Mannheim, Germany}

\title{Regression based thresholds in principal loading analysis}

%Threshold values in PLA from an OLS regression point of view
%Threshold values in PLA according to an OLS regression view
%Threshold choice in principal loading analysis based on an OLS regression approach
%Threshold choice in principal loading analysis based on an OLS regression view
%Threshold choice in principal loading analysis according to an OLS regression sense

\begin{abstract}
Principal loading analysis is a dimension reduction method that discards variables which have only a small distorting effect on the covariance matrix. As a special case, principal loading analysis discards variables that are not correlated with the remaining ones. In multivariate linear regression on the other hand, predictors that are neither correlated with both the remaining predictors nor with the dependent variables have a regression coefficients equal to zero. Hence, if the goal is to select a number of predictors, variables that do not correlate are discarded as it is also done in principal loading analysis. That both methods select the same variables occurs not only for the special case of zero correlation however. We contribute conditions under which both methods share the same variable selection. Further, we extend those conditions to provide a choice for the threshold in principal loading analysis which only follows recommendations based on simulation results so far. \\[12pt]
\keywords{Dimensionality Reduction, Matrix Perturbation Theory, Ordinary Least Squares Regression, Multivariate Linear Regression, Principal Loading Analysis, Principal Component Analysis}
\\[12pt]
\mathsubjclassSingle{2010 Mathematics Subject Classification}{
62J05,  	%Linear regression; mixed models
47H14,  	%Perturbations of nonlinear operators
62H25, %Factor analysis and principal components; correspondence analysis
15A18, %Eigenvalues, singular values, and eigenvectors
}
\end{abstract}
\end{frontmatter}

%% main text
\section{Introduction}
\label{s:Introduction}
\noindent
Principal loading analysis (PLA) is a tool developed by \citet{BD21} to reduce dimensions. The method chooses a subset of observed variables by discarding the other variables based on the impact of the eigenvectors on the covariance matrix. \citet{BA21} extended the method to correlation-based PLA in order to provide a scale invariant approach.

By nature, some parts of PLA correspond with principal component analysis (PCA) which is a popular dimension reduction technique first formulated by \citet{PE01} and \citet{HO33}. Both are based on the eigendecomposition of the covariance (correlation) matrix and take into account the magnitude of dispersion of each eigenvector. Despite the intersections however, the outcome is different because PCA yields a reduced set of variables, named principal components (PCs), by transforming the original variables.

PCA has been also extended to ordinary least squares (OLS) regression contexts by \citet{HO57} and \citet{KE57}. Instead of regressing on the original variables, in principal component regression the PCs, which are orthogonal by construction, are used as predictors to address multicollinearity. There also exist approaches to reduce the number of PCs based on principal component regression (see \citet{HO76} and \citet{MW77} among others).

However, even if reducing the number of PCs yields yet another reduced set of all original variables and not a subset of the original variables, it is natural to examine an extension to linear regression contexts also for PLA. We provide such an extension to multivariate linear regression (MLR). Linear regression is then a special case of our results with reduced dimension. Two things have to be taken into account: firstly, MLR can be interpreted as a dimensionality reduction tool itself when predictors are eliminated if the corresponding regression coefficients are not significantly different from zero. Hence, a subset of the original variables is selected by discarding the non-significant variables, which can be done by backwards elimination. Secondly, MLR considers predictors $\Xb_\bK$ that are regressed on independent variables $\Xb_\bD$, while PLA considers cross-sectional data. To address this difference, PLA is conducted on $\Xb \equiv (\Xb_\bK,\Xb_\bD)$ in this work to take into account both the predictors as well as the independent variables. The idea to perform the method not only on $\Xb_\bD$ but rather on $\Xb$ is not new and has been used for latent root regression, a variation of PC regression, which has been formulated by \citet{HA73} and \citet{WGM74}, however, not in a multivariate regression context.

In this work, we check for intersections of MLR and PLA in a dimensionality reduction-sense. That is, we examine the conditions under which the same variables are selected hence the same variables are discarded by MLR and PLA. We provide necessary as well as necessary and sufficient conditions for those intersections. \citet{AN63}, \citet{KN93} and \citet{NW90} provide convergence rates for the sample noise of covariance matrices which are useful for our results because PLA and MLR only share the same selection of variables in the sample case and not in the population case except for the trivial case of no perturbation of the covariance matrix. We interpret sample values as their perturbed population counterparts. Hence, our results are based on matrix perturbation theory which originally comes from algebraic fields such as in \citet{IN09} and \citet{SS90}. The contributed bounds  for the cut-off value are based on the Davis-Kahan theorem (\citet{DK70} and \citet{YW15}) as it was done in the work of \citet{BD21}.

This article is organized as follows: Section~\ref{s:setup} provides notation needed for the remainder of this work. In Section~\ref{s:methodology}, we recap PLA and introduce the concept of our approach for comparison. An intuition for our conclusions is provided by Lemma~\ref{th:intuition} in Section~\ref{s:PLAtoOLS} which is the base for all remaining results. Further, we contribute bounds for the perturbation matrices under which PLA and MLR regression select the same variables and under which they yield different results. Same is done for the threshold value in PLA in Section~\ref{s:Threshold}. We give an example in Section~\ref{s:example} and take a resume as well as suggest extensions in Section~\ref{s:conclusion}.

\section{Setup}
\label{s:setup}

We first state some notation and assumptions. Generally, we use index sets to assign submatrices which is extensively needed in this work. For any set $\bK = \{k_1,\ldots,k_K \} \subset \N$ with $k_1<\ldots<k_K$, $\Ab_{[\{r\},\bK]}$ denotes the $r$th row of columns $k_1,\ldots,k_K$ for a matrix $\Ab$. For columns, $\Ab_{[\bK,\{c\}]}$ follows the same notation. Further, we denote the elements of any column vector $\vb$ by $\vb = ( v^{(1)} , \ldots ,  v^{(M)} )$.

We consider $\xb = (\xb_1,\ldots,\xb_M) \equiv  (\xb_\bK , \xb_\bD)$ $\in\R^{N\times M}$ to be an independent and identically distributed (IID) sample. $\xb$ contains $n\in\{1,\ldots,N\}$ with $N>M$ observations of a random vector $\Xb = (X_1,\ldots,X_M) \equiv (\Xb_\bK, \Xb_\bD)$ with $\Xb_\bK^\top \in \R^K$ and $\Xb_\bD^\top \in \R^D$ such that $K + D = M$. The covariance matrix of $\Xb$ is given by $\Sigmab = (\sigma_{i,j})$ for defined indices $i,j\in\{1,\ldots,M\}$ throughout this work. Further, $\Sigmab_1$ and $\Sigmab_2$ are the covariance matrices of $\Xb_\bK$ and $\Xb_\bD$ respectively such that
\begin{equation} \label{eq:covblock}
\Sigmab \equiv \begin{pmatrix}
\Sigmab_1 & \Sigmab_{12} \\ \Sigmab_{12}^\top  & \Sigmab_2\end{pmatrix}\;.
\end{equation}
 Analogously to $\xb$, we introduce the IID sample $\xbtilde = (\xbtilde_1,\ldots,\xbtilde_M) \equiv  (\xbtilde_\bK , \xbtilde_\bD)$ $\in\R^{N\times M}$ drawn from a random vector $\Xbtilde = (\Xtilde_1,\ldots,\Xtilde_M) \equiv (\Xbtilde_\bK, \Xbtilde_\bD)$ with covariance matrix $\Sigmabtilde$. The intuition is that $\Xb$ and $\Xbtilde$, hence also $\xb$ and $\xbtilde$, behave similarly, however, with the distinction that $\Sigmabtilde \equiv \Sigmab + \Epsilonb$ is slightly perturbed by a sparse matrix $\Epsilonb = (\varepsilon_{i,j})$. $\Epsilonb$ is a technical construction that contains small components extracted from $\Sigmab$ hence $\varepsilon_{i,j} \neq 0 \Rightarrow \sigma_{i,j} = 0$. The sample counterparts of $\Sigmab$ and $\Sigmabtilde$ are given by
\begin{equation*}
\Sigmabhat \equiv \Sigmab + \Etab
\end{equation*}
and
\begin{equation*}
\Sigmabtildehat  \equiv \Sigmabtilde +  \Etabtilde \equiv \Sigmab + \Epsilonb + \Etabtilde 
\end{equation*}
respectively. $\Etab$ and $\Etabtilde$ are perturbations in the form of random noise matrices. The noise is due to a finite number of sample observations. $\Sigmabhat$, $\Sigmabtildehat$, and all of their compound matrices follow an analogous block-structure as $\Sigmab$ given in (\ref{eq:covblock}).

Further, we consider the MLR models $\xb_\bK = \xb_\bD \betab + \ub$ and $\xbtilde_\bK = \xbtilde_\bD \betabtilde + \ubtilde$ with $\ub,\ubtilde\in\R^{N\times K}$, $\betab = (\betab_1,\ldots,\betab_K) \in \R^{D\times K}$ and $\betabtilde = (\betabtilde_1,\ldots,\betabtilde_K)\in \R^{D\times K}$ where $\betab_k,\betabtilde_k \in \R^D$ for $k\in\{1,\ldots,K\}$. We denote the estimators for $\betab$ and $\betabtilde$ by $\betabhat$ and $\betabtildehat$ respectively.

We use the $\mathcal{L}_1$ and $\mathcal{L}_2$ vector norms $\Vert\vb\Vert_1 \equiv \sum_i |v^{(i)}|$ and $\Vert\vb\Vert_2 \equiv \left(\sum_i (v^{(i)})^2\right)^{1/2}$, and the Frobenius norm $\Vert\Ab\Vert_F$ as a subordinate extension for matrices. Further, the trace of a matrix $\Ab$ is given by $\tr(\Ab)$. $\bigO(\cdot)$ denotes stochastic boundedness and we use the symbol "$\ind$" as a shortcut that variables have zero correlation and therefore not necessarily for independence. Further, $\bm{0}$ denotes a matrix containing zeros and its dimension its always self-explanatory.

Variables or blocks of variables are linked to eigenvectors in PLA. Hence, for practical purposes we define $\bD \subset \{1,\ldots,M\}$ with $\card(\bD) = D$ such that $1\leq D  < M$. $\bD$ will contain the indices of a variable-block $\{\Xtilde_d\} \equiv \{\Xtilde_d\}_{d\in\bD}$ we discard according to PLA, and we introduce the quasi-complement $\bK \equiv  \{1,\ldots,M\} \backslash \bD$ with $\card(\bK) = K$ to label the remaining variables that are kept. In an analogue manner, $\bDelta$ with $\card(\bDelta) = D$ will be used to index eigenvectors with respective eigenvalues linked to the $\{\Xtilde_d\}$. $\bD\times\bDelta \equiv \{(d,\delta) : d\in\bD, \; \delta\in\bDelta\}$ denotes the Cartesian product and therefore contains all possible combinations of paired indices $d$ linked to $\delta$.  For convenience purposes, we consider only a single block $\{\Xtilde_d\}$ to be discarded throughout this work except for the example in Section~\ref{s:example}. However, our results remain true for any number of blocks.

Given the notation above, the following assumptions are made:

\begin{assumption}{1}{}\label{ass:independence}
$\Xb_\bK \ind \Xb_\bD$.
\end{assumption}
Zero correlation among the random variables is a technical construction that allows to extract the correlations among $(\Xtilde_1,\ldots,\Xtilde_{M})$ using the perturbation $\Epsilonb$ of the covariance matrix. Hence, we do not assume that there is zero correlation among the random variables $(\Xtilde_1,\ldots,\Xtilde_{M})$.

\begin{assumption}{2}{}\label{ass:meancentred}
$\E ( \Xb )= \E( \Xbtilde )  = \bm{0}$ and the samples $\xb$ as well as $\xbtilde$ have mean zero.
\end{assumption}
We assume variables to have mean zero for convenience purposes to simplify the estimators of the covariance matrices to $\Sigmabhat = (N-1)^{-1}\xb^\top\xb$ and $\Sigmabtildehat = (N-1)^{-1}\xbtilde^\top\xbtilde$ respectively.

\begin{assumption}{3}{}\label{ass:OLS}
The sample regression coefficients $\betabhat$ and $\betabtildehat$ of the linear models $\Xb_\bK = \Xb_\bD \betab + \ub$ and $\Xbtilde_\bK = \Xbtilde_\bD \betabtilde + \ubtilde$ are estimated according to $\betabhat = \Sigmabhat_2^{-1} \Sigmabhat_{12}^\top$ and $\betabtildehat =\Sigmabtildehat_2^{-1} \Sigmabtildehat_{12}^\top$ respectively. The error terms $\ub$ and $\ubtilde$ are mean zero with covariance matrices $\Sigmab_{\ub}$ and $\Sigmabtilde_{\ub}$.
\end{assumption}

It is assumed that the regression coefficients can be estimated using OLS. The respective sample linear models are $\xb_\bK = \xb_\bD \betabhat + \ub$ and $\xbtilde_\bK = \xbtilde_\bD \betabtildehat + \ubtilde$ with  $\betabhat = \Sigmabhat_2^{-1} \Sigmabhat_{12}^\top \equiv (\xb_\bD^\top\xb_\bD)^{-1} \xb_\bD^\top \xb_\bK$ and $\betabtildehat = \Sigmabtildehat_2^{-1} \Sigmabtildehat_{12}^\top \equiv (\xbtilde_\bD^\top\xbtilde_\bD)^{-1} \xbtilde_\bD^\top \xbtilde_\bK$ due to Assumption~\ref{ass:meancentred}. If the regression coefficient of a variable is not-significantly different from zero, we say that the respective variable is a non-significant predictor and that the variable can be discarded according to MLR. We implicitly assume that the fourth-order moments $ \E( \Xb^\top\Xb \otimes \Xb^\top\Xb) < \infty$ and $ \E( \Xbtilde^\top\Xbtilde \otimes \Xbtilde^\top\Xbtilde) < \infty$ are finite because OLS estimation is sensitive to outliers. This further allows us to obtain the rate of convergence for the sample covariance matrices by applying the Lindeberg-L\'evy central limit theorem.

\begin{assumption}{4}{}\label{ass:eigenvalues}
The eigenvalues $\lambda_1 > \ldots > \lambda_M > 0$ of $\Sigmab$ and the eigenvalues $\lambdatildehat_1 > \ldots > \lambdatildehat_M > 0$ of $\Sigmabtildehat$ are arranged in descending order.
\end{assumption}

This assumption is in line with the work of \citet{BD21} where the asymptotics of eigenvalues were considered. The reason behind this is, that asymptotics for eigenvalues with a multiplicity greater than one are problematic to obtain as pointed out by \citet{DPR82}. In this work, we are not concerned with the limits of eigenvalues. However, we make usage of the Davis-Kahan Theorem and hence divide by the difference of eigenvalues. Assumption~\ref{ass:eigenvalues} therefore is rather strict, but by loosing it one has to assume that the difference of the respective eigenvalues in Section~\ref{s:Threshold} is unequal to zero. Further, we assume that the population covariance matrix $\Sigmab$ and the sample covariance matrix $\Sigmabtildehat$ of $\xbtilde$ are positive definite and hence invertible.
\section{Methodology}
\label{s:methodology}

In this section, we recap PLA and state our general approach to derive the conditions under which PLA and MLR share the same selection of variables. The following summary of PLA is based on the work of \citet{BD21} and we refer to this for a more elaborate explanation.

PLA is a tool for dimension reduction where a subset of existing variables is selected while the other variables are discarded. The intuition is that blocks of variables are discarded which distort the covariance matrix only slightly.  This is done by checking if all $K$ elements $k\in\bK$ of $D$ eigenvectors $\{\vtilde_{\delta}^{(k)}\}_{\delta\in\bDelta}$ of the covariance matrix $\Sigmabtilde$ are smaller in absolute terms than a certain cut-off value $\tau$. That is to say, we check if $|\vtilde_{\delta}^{(k)}| \leq \tau$ for all $(k,\delta)\in\bK\times\bDelta$. If this is the case, then most distortion of the variables $\{\Xtilde_d\}$ is geometrically represented by the eigenvectors $\{\vbtilde_{\delta}\}_{\delta\in\bDelta}$. The explained variance of the $\{\Xtilde_d\}$ can then be evaluated either by $N^{-1}\sum_{d\in\bD} \sigmatilde_{d,d}$ or approximated by $(  \sum_{m\in\bP} \lambdatilde_m )^{-1} ( \sum_{\delta\in\bDelta} \lambdatilde_{\delta} ) $ where $\lambdatilde_i$ is the eigenvalue corresponding to the eigenvector $\vbtilde_i$. If the explained variance is sufficiently small for the underlying purpose of application, the variables $\{\Xtilde_d\}$ are discarded. In practice, the population eigenvectors and eigenvalues are replaced by their sample counterparts.

For now, we consider the extreme case that a variable is discarded according to PLA using a threshold $\tau = 0$. According to this, the eigenvectors need to contain zeros which reflects zero correlation among variables. In MLR on the other hand, a predictor that is neither correlated with the remaining predictors nor with the dependent variables has regression coefficients equal to zero. Hence, if the goal is to select a number of predictors, the variable that is not correlated is discarded. Therefore, there is a natural link between PLA and MLR based on the correlation structure and we cover those underlying structures that result in the same selections of variables, both according to PLA and MLR.

We do not take the explained variance of blocks into account despite its importance for PLA. The reason is that it is only relevant if predictors lie in the same block as the independent variable to select them in a PLA sense.

\section{Perturbation Bounds}
\label{s:PLAtoOLS}

In this section, we derive conditions for the perturbation matrices that have to hold in order to obtain the same selected variables by PLA and MLR. To be more precise: we provide conditions under which discarded variables according to PLA are also considered to be non-significant predictors in MLR. Here, a non-significant predictor means that the regression coefficients of the respective variables are not significantly different from zero.

The following theorem provides an intuition of the relationship between the perturbation matrices:

\begin{lemma}\label{th:intuition}
Let the variables $\{\Xtilde_d\}$ be discarded according to PLA and let $\{X_d\}$ be non-significant predictors. $\Xtilde_\dstar \in \{\Xtilde_d\}$ is considered to be a non-significant predictor if for all $k\in\bK$
\begin{equation*}
1\geq \dfrac{ (\Sigmab_2^{-1})_{[\{\dstar\},\bD]} ( \Epsilonb_{12}+ \Etabtilde_{12} )_{[\bD,\{k\}]}  }{ (\Sigmab_2^{-1})_{[\{\dstar\},\bD]}   {\Etab_{12}}_{[\bD,\{k\}]}  } \geq -1 \;  .
\end{equation*}
\end{lemma}

The intuition is that the perturbation together with the noise of $\Sigmabtilde$ shall not billow more than pure noise itself. Simply put, this can either happen in the trivial case when $\Epsilonb = \bm{0}$, or when $\Epsilonb$ and $\Etabtilde$ fluctuate towards different directions hence pull themselves back below $\Etab$. However, strictly spoken the perturbations are sandwiched hence we have to take multiplication with $(\Sigmab_2^{-1})_{[\{\dstar\},\bD]} $ into account. Nonetheless, Lemma~\ref{th:intuition} provides an intuition regarding the perturbations.\\
More detailed results regarding the perturbation and noise are given in Theorem~\ref{l:bounds} and in Corollary~\ref{c:bounds}. The results provide necessary as well as necessary and sufficient conditions under which PLA and MLR share the same selection of variables.

\begin{theorem}\label{l:bounds}
Let the variables $\{\Xtilde_d\}$ be discarded according to PLA and let $\{X_d\}$ be non-significant predictors. $\Xtilde_\dstar \in \{\Xtilde_d\}$ is considered to be a non-significant predictor if for all $k\in\bK$
\begin{equation*}
 \Vert ( \Epsilonb_{12}^\top+ \Etabtilde_{12}^\top )_{[\bD,\{k\}]}  \Vert_2 \leq |   \dfrac{ (\Sigmab_2^{-1})_{[\{\dstar\},\bD]}  }{ \Vert (\Sigmab_2^{-1})_{[\{\dstar\},\bD]} \Vert_2 }  (\Etab_{12}^\top)_{[\bD,\{k\}]}   |    \;.
\end{equation*}
\end{theorem}

If $\Xtilde_\dstar$ is discarded according to PLA, the variable was weakly correlated with all $\{\Xtilde_k\}_{k\in\bK}$. Further, the smaller the partial correlation of $\Xtilde_\dstar$ within $\{\Xtilde_d\}$, the more likely $\Xtilde_\dstar$ is also considered to be a non-significant predictor. In practice, $(\Sigmab_2^{-1})_{[\{\dstar\},\bD]}$ can be estimated by $(\Sigmabhat_2^{-1})_{[\{\dstar\},\bD]} $. We can further bound the perturbations of $\Sigmabtilde$ by noise of $\Sigmab$ only, which reflects the intuition discussed in Lemma~\ref{th:intuition}.

\begin{corollary}\label{c:bounds}
Let the variables $\{\Xtilde_d\}$ be discarded according to PLA and let $\{X_d\}$ be non-significant predictors. Necessary conditions that the $\{\Xtilde_d\}$ are considered to be non-significant predictors are
\begin{equation*}
 \Vert ( \Epsilonb_{12}^\top + \Etabtilde_{12}^\top)_{[\bD,\{k\}]}  \Vert_2 \in [ 0 , \Vert (\Etab_{12}^\top)_{[\bD,\{k\}]} \Vert_2 ]   \text{ for all }k\in\bK \; .
\end{equation*}
\end{corollary}

Corollary~\ref{c:bounds} reflects what we discussed above, that is that the perturbation together with the noise of $\Sigmabtilde$ shall not billow more than pure noise itself. Further, it is indicated that there are conditions under which PLA and MLR select the same variables in the sample case, however not necessarily in the population case. This is due to the perturbation $\Epsilonb$ as stated in the following result.

\begin{corollary}\label{c:convergencerate}
Let the variables $\{\Xtilde_d\}$ be discarded according to PLA and let $\{X_d\}$ be non-significant predictors. Further, let  $\Xtilde_\dstar \in \{\Xtilde_d\}$ considered to be a non-significant predictor. Then $\Epsilonb = \bigO(N^{-1/2})$.
\end{corollary}

\citet{BD21} provided convergence rates for the noise present in all sample parts needed for PLA. However, a bound for the perturbation $\Epsilonb$ was not given. Looking from an MLR perspective onto PLA, $\Epsilonb$ has to be bounded by $\bigO(N^{-1/2})$. For the trivial case that $\Epsilonb = \bm{0}$, MLR and PLA share the same selection of variables also for the population case. When perturbation is present however, that is to say when the variables to be discarded according to PLA are correlated with the remaining variables, both methods can only select the same variables in the finite sample case.
\section{Implications for the Threshold Value}
\label{s:Threshold}

Choosing an appropriate threshold $\tau$ is crucial for PLA. \citet{BD21} and \citet{BA21} provide feasible thresholds for PLA based on simulations. In this section, we translate our results from Section~\ref{s:PLAtoOLS} to the choice of such a cut-off value. We provide bounds for the threshold $\tau$ that ensure discarding according to PLA as well as testing the respective regression coefficient to be not significantly different from zero on a $(1-\alpha)\cdot 100\%$ level hence ensure discarding in a multivariate regression sense. Therefore, the threshold depends on the significance level $\alpha$. Further, we provide bounds for the perturbations $\Epsilonb + \Etabtilde$ in this section.

The procedure to test if all regression coefficients are significantly different from zero is well known. Commonly, the tests based on Roy's largest root, Wilks' Lambda, Hotelling-Lawley trace and Pillai-Bartlett trace are used. In a multivariate regression sense, all variables $\{\Xtilde_d\}$ are discarded if we fail to reject the null of $\betab = \bm{0}$.

\begin{theorem}\label{l:tau}
Let $\tau$ be the cut-off value for PLA and $\lambda_1 > \ldots > \lambda_{M}$ the eigenvalues of $\Sigmab$ with $\lambda_0 \equiv \infty$ and $\lambda_{M+1} \equiv - \infty$. $\{\Xtilde_d\}$ are discarded both by MLR on a $(1-\alpha)\cdot 100\%$ level and by PLA if
\begin{equation*}
\dfrac{2^{2/3}\Vert \  \Epsilonb+ \Etabtilde  \Vert_F }
{\min_{\delta\in\bDelta} \{\lambda_{\delta-1} - \lambda_\delta, \lambda_\delta - \lambda_{\delta+1}\}} 
\leq \tau \leq 
 \dfrac{2^{2/3}   \Vert \Sigmabtildehat_1^{-1} \Vert_F^{-1/2} \Vert \Sigmabtildehat_2^{-1} \Vert_F^{-1/2}  \dfrac{S D_1}{F_\alpha^{-1}D_2 + D_1}  }
 {\min_{\delta\in\bDelta} \{\lambda_{\delta-1} - \lambda_\delta, \lambda_\delta - \lambda_{\delta+1}\}} 
\end{equation*}
where $F_\alpha$ is the quantile of the $F$-distribution with degrees of freedom $D_1 \equiv D  K$ and $D_2 \equiv S (N-D-K)$ with $S\equiv \min\{D,K\}$.
\end{theorem}

Choosing $\tau_\alpha \equiv \dfrac{2^{2/3}  F_\alpha \Vert \Sigmabtildehat_1^{-1} \Vert_F^{-1/2} \Vert \Sigmabtildehat_2^{-1} \Vert_F^{-1/2}  S D_1 D_2^{-1} }{\min_{\delta\in\bDelta} \{\lambda_{\delta-1} - \lambda_\delta, \lambda_\delta - \lambda_{\delta+1}\}}  $ is sufficient, not necessary and sufficient. This is because the terms are bound from above in the derivations.

Since $\Sigmabtildehat \equiv \Sigmab + \Epsilonb + \Etabtilde$, it is difficult to obtain bounds for $\Epsilonb$ and $\Etabtilde$ having only access to $\Sigmabtildehat$ in practice. However, if the $\{\Xtilde_d\}$ are not discarded according to multivariate regression, we can find a lower bound for $\Epsilonb+ \Etabtilde$.

\begin{corollary}\label{c:perturbationbounds}
If the $\{\Xtilde_d\}$ are not discarded by MLR on a $(1-\alpha)\cdot 100\%$ level, it holds hat
\begin{equation*}
\Vert \  \Epsilonb+ \Etabtilde  \Vert_F > 2^{2/3}   \Vert \Sigmabtildehat_1^{-1} \Vert_F^{-1/2} \Vert \Sigmabtildehat_2^{-1} \Vert_F^{-1/2}  \dfrac{S D_1}{F_\alpha^{-1}D_2 + D_1}
\end{equation*}
where $F_\alpha$ is the quantile of the $F$-distribution with degrees of freedom $D_1 \equiv D  K$ and $D_2 \equiv S (N-D-K)$ with $S\equiv \min\{D,K\}$.
\end{corollary}

For completion, we further give necessary conditions in line with the intuition from Section~\ref{s:PLAtoOLS}.

\begin{corollary}\label{c:tau}
Let $\tau$ be the cut-off value for PLA and $\lambda_1 > \ldots > \lambda_{M}$ the eigenvalues of $\Sigmab$ with $\lambda_0 \equiv \infty$ and $\lambda_{M+1} \equiv - \infty$. Let further $\{X_d\}$ be non-significant predictors. A necessary condition to discard $\{\Xtilde_d\}$ both by PLA and by MLR is given by
\begin{equation*}
\dfrac{2^{2/3}\Vert \  \Epsilonb+ \Etabtilde \Vert_F }{\min_{\delta\in\bDelta} \{\lambda_{\delta-1} - \lambda_\delta, \lambda_\delta - \lambda_{\delta+1}\} }  \leq \tau \leq \dfrac{2^{2/3}\Vert \   (\Etab_{12}^\top)_{[\bD,\bK]}   \Vert_2  }{\min_{\delta\in\bDelta} \{\lambda_{\delta-1} - \lambda_\delta, \lambda_\delta - \lambda_{\delta+1}\} }  \;.
\end{equation*}
\end{corollary}

Note that the case when we focus only on a single variable which belongs to a block $\{\Xtilde_d\}$ in PLA-sense instead of considering the whole block itself, the maximization operator in Theorem~\ref{l:tau} simplifies. The reason is that the link between the variable and the eigenvector with respective eigenvalue is evident in this case. For a block of variables it is not obvious which variable is linked to which eigenvectors and hence all eigengaps $ \{\lambda_{\delta-1} - \lambda_\delta, \lambda_\delta - \lambda_{\delta+1}\}$ for $\delta\in\bDelta$ have to be considered.

\begin{corollary}\label{c:tau2}
Let $\tau$ be the cut-off value for PLA and $\lambda_1 > \ldots > \lambda_{M}$ the eigenvalues of $\Sigmab$ with $\lambda_0 \equiv \infty$ and $\lambda_{M+1} \equiv - \infty$. $\Xtilde_\dstar \in \{\Xtilde_d\}$ is discarded both by MLR on a $(1-\alpha)\cdot 100\%$ level and by PLA if
\begin{equation*}
\dfrac{2^{2/3}\Vert \  \Epsilonb+ \Etabtilde  \Vert_F }
{\min\{\lambda_{\deltastar-1} - \lambda_\deltastar, \lambda_\deltastar - \lambda_{\deltastar+1}\}} 
\leq \tau \leq 
 \dfrac{2^{2/3}   \Vert \Sigmabtildehat_1^{-1} \Vert_F^{-1/2} \Vert \Sigmabtildehat_2^{-1} \Vert_F^{-1/2}  \dfrac{S D_1}{F_\alpha^{-1}D_2 + D_1} }
 {\min\{\lambda_{\deltastar-1} - \lambda_\deltastar, \lambda_\deltastar - \lambda_{\deltastar+1}\}} 
\end{equation*}
where $F_\alpha$ is the quantile of the $F$-distribution with degrees of freedom $D_1 \equiv D  K$ and $D_2 \equiv S (N-D-K)$ with $S\equiv \min\{D,K\}$, and $\vbtilde_\deltastar$ is the eigenvector linked to $\Xtilde_\dstar$.
\end{corollary}

Same holds also for Corollary~\ref{c:tau} which also simplifies when considering only a single variable.

\begin{corollary}\label{c:tau3}
Let $\tau$ be the cut-off value for PLA and $\lambda_1 > \ldots > \lambda_{M}$ the eigenvalues of $\Sigmab$ with $\lambda_0 \equiv \infty$ and $\lambda_{M+1} \equiv - \infty$.  If $\{X_d\}$ are non-significant predictors, then a necessary condition to discard $\Xtilde_{\dstar}\in\{\Xtilde_d\}$ both by PLA and by MLR is given by
\begin{equation*}
\dfrac{2^{2/3}\Vert \  \Epsilonb+ \Etabtilde \Vert_F }{\min \{\lambda_{\deltastar-1} - \lambda_\deltastar, \lambda_\deltastar - \lambda_{\deltastar+1}\} }  \leq \tau \leq \dfrac{2^{2/3}\Vert \   (\Etab_{12}^\top)_{[\bD,\bK]}   \Vert_2  }{\min \{\lambda_{\deltastar-1} - \lambda_\deltastar, \lambda_\deltastar - \lambda_{\deltastar+1}\} }  
\end{equation*}
where $\vbtilde_\deltastar$ is the eigenvector linked to $\Xtilde_\dstar$.
\end{corollary}

In practice, the eigenvalues $\lambda_\delta$ can be estimated by $\lambdatildehat_\delta$. \citet{BD21} outlined that estimation is biased by $\Vert \Epsilonb \Vert_F$ which is small, however, due to the sparseness of $\Epsilonb$.

\section{Example}
\label{s:example}
We provide an example how multivariate regression can support the threshold choice in PLA. Firstly, we conduct PLA with a threshold $\tau$ to detect the underlying blocks. Afterwards, the threshold is adjusted to $\tau_\alpha$ according to a MLR-sense. The blocks can then either be discarded by PLA with threshold $\tau_\alpha$ as well, which corresponds to the case that both methods choose the same variables, or are not discarded by PLA with threshold $\tau_\alpha$. The latter one corresponds to the case where variables are discarded by PLA with threshold $\tau$ but not according to MLR. Further, we give an example where a variable is discarded by MLR, however, not by PLA. All values are rounded to four decimal places.

We assume a sample to be drawn from the random vector $(\Xtilde_1, \Xtilde_2, \Xtilde_3, \Xtilde_4, \Xtilde_5)$ and we generated a respective sample $\xbtilde$ of $N=100$ observations. The sample covariance matrix $\Sigmabtildehat$ can be found in \ref{a:compl} while the eigenvectors $\vbtildehat_m$ and eigenvalues $\lambdatildehat_m$ with $m\in\{1,\ldots,5\}$ are in Table~\ref{table:example}.

\begin{table}[h!]
\centering
\caption{Eigenvectors $\vbtildehat_m$ and eigenvalues $\lambdatildehat_m$ with $m\in\{1,\ldots,5\}$ of the sample covariance matrix $\Sigmabtildehat$. All values above $\tau = 0.3$ are highlighted in bold.}
\label{table:example}
\begin{tabular}{lrrrrr}
\toprule
$\vbtildehat_m$ & $\vbtildehat_1$ & $\vbtildehat_2$ &$\vbtildehat_3$ & $\vbtildehat_4$ & $\vbtildehat_5$ \\
  \cmidrule(l{3pt}r{3pt}){1-1} \cmidrule(l{3pt}r{3pt}){2-6}
$\Xtilde_1$ & -0.1143 &  \textbf{0.9578} & -0.0782 &  0.2435 & -0.0642 \\
$\Xtilde_2$ & -0.1453 & -0.0777 & -0.0310 & -0.0326 & \textbf{-0.9853} \\
$\Xtilde_3$ &  0.1731 & -0.1106 &  \textbf{0.6664} &  \textbf{0.7141} & -0.0614 \\
$\Xtilde_4$ & \textbf{ 0.5797} &  0.2511 &  \textbf{0.5036} & \textbf{-0.5804} & -0.1019 \\
$\Xtilde_5$ & \textbf{-0.7744} &  0.0364 &  \textbf{0.5433} &  \textbf{-0.3047} &  0.1043 \\
  \cmidrule(l{3pt}r{3pt}){1-1} \cmidrule(l{3pt}r{3pt}){2-6}
$\lambdatildehat_m$ & 10.2370 &  8.9338 &  6.0507 &  4.9167 &  1.1202 \\
\toprule
\end{tabular}
\end{table}

We choose $\tau = 0.3$ as recommended by \citet{BD21}, and $\alpha = 0.05$ for our analysis. PLA detects three Blocks: $\Xtilde_1$ ($1\times1$), $\Xtilde_2$ ($1\times1$) and $(\Xtilde_3,\Xtilde_4,\Xtilde_5)$ ($3\times3$) as can be seen in Table~\ref{table:example}. We now want to bear the PLA results from a MLR point of view and check the two $1\times1$ blocks $\Xtilde_1$ and $\Xtilde_2$ according to Corollary~\ref{c:tau2}.

Starting with $\Xtilde_1$, we consider the MLR model $(\xbtilde_2,\xbtilde_3,\xbtilde_4,\xbtilde_5) = \xbtilde_1 \betabtilde + \ubtilde$. Hence, $D_1 \equiv 1 \cdot 4 =4$, $D_2 \equiv 1 \cdot (100 - 1 - 4) = 95$, and $F_\alpha =  2.4675$. The block is reflected by the second eigenvector and therefore we have that $\delta = 2$. Since $\Vert \Sigmabtildehat_1^{-1} \Vert_F^{-1/2} = 0.9597$, $\Vert \Sigmabtildehat_2^{-1} \Vert_F^{-1/2} = 0.4299$ (\ref{a:compl}), and $\min\{  \lambdatildehat_1 - \lambdatildehat_2, \lambdatildehat_2 - \lambdatildehat_3\}  = 1.3032$, we can conclude that a sufficient threshold to discard $\Xtilde_1$ by MLR is given by $\hat{\tau}_{\alpha,1} =  0.2753$. From Table~\ref{table:example} we can see that PLA detects the block containing $\Xtilde_1$ also for $\hat{\tau}_{\alpha,1}$ making it a feasible threshold for both, PLA and MLR. Therefore, the result we have gotten from PLA in the first place is also supported by regression analysis.

On the other hand, for the block consisting of $\Xtilde_2$ we consider the MLR model $(\xbtilde_1,\xbtilde_3,\xbtilde_4,\xbtilde_5) = \xbtilde_2 \betabtilde + \ubtilde$ with $\Vert \Sigmabtildehat_1^{-1} \Vert_F^{-1/2} = 0.5568$, $\Vert \Sigmabtildehat_2^{-1} \Vert_F^{-1/2} = 0.9326$ (\ref{a:compl}), and $\min\{  \lambdatildehat_4 - \lambdatildehat_5, \lambdatildehat_5 - \lambdatildehat_6\}  = \lambdatildehat_4 - \lambdatildehat_5 = 3.7965$. $\delta = 5$ because $\Xtilde_2$ is reflected by the fifth eigenvector. We can conclude that a sufficient threshold to discard $\Xtilde_2$ by MLR is given by $\hat{\tau}_{\alpha,2} =  0.0751$. Using this threshold for PLA however, $\Xtilde_2$ is not detected as a block. Hence, this block is not supported by MLR and we can conclude that $\Vert \  \Epsilonb+ \Etabtilde  \Vert_F >  \Vert \Sigmabtildehat_1^{-1} \Vert_F^{-1/2} \Vert \Sigmabtildehat_2^{-1} \Vert_F^{-1/2}  \dfrac{S D_1}{D_2 \cdot F_\alpha^{-1} +D_1} \approx 0.1795$ according to Corollary~\ref{c:perturbationbounds}

Considering the regression model  $(\xbtilde_1,\xbtilde_2,\xbtilde_4,\xbtilde_5) = \xbtilde_3 \betabtilde + \ubtilde$, the estimated regression coefficients $\betabtildehat$ are not significantly different from $\bm{0}$ on a 5\% level. In fact, the $p$-values for the approximated $F$-tests based on Roy's largest root, Wilks' Lambda, Hotelling-Lawley trace and Pillai-Bartlett trace all equal 0.4685 with the respective value of the test statistic being $0.8978$ (\ref{a:compl}). Hence, $\Xtilde_3$ is discarded in a MLR sense. However, $\Xtilde_3$ is not considered to be a block in PLA-sense which means that $  \Vert \Sigmabtildehat_1^{-1} \Vert_F^{-1/2} \Vert \Sigmabtildehat_2^{-1} \Vert_F^{-1/2}  \dfrac{S D_1}{D_2 \cdot F_\alpha^{-1} +D_1} \approx 0.2259 >  \Vert \  \Epsilonb+ \Etabtilde  \Vert_F > 2^{-2/3} \tau \cdot \min\{  \lambdatildehat_{\delta-1} - \lambdatildehat_\delta, \lambdatildehat_\delta - \lambdatildehat_{\delta+1}\}^{-1}$ (\ref{a:compl}). Since $\Xtilde_3$ is not detected to be a $1\times 1$ block according to PLA, it is difficult to choose the eigenvector that represents $\Xtilde_3$ in order to calculate the respective eigengaps. However, it is reasonable to assume the vector to be either $\vbtildehat_3$ or $\vbtildehat_4$ because  $\Xtilde_3$ obtains the largest values in them. Since $\min_{\delta=3} \{  \lambdatildehat_{\delta-1} - \lambdatildehat_\delta, \lambdatildehat_\delta - \lambdatildehat_{\delta+1}\} = \lambdatildehat_3 - \lambdatildehat_{4} = 1.1340$ and $\min_{\delta=4} \{  \lambdatildehat_{\delta-1} - \lambdatildehat_\delta, \lambdatildehat_\delta - \lambdatildehat_{\delta+1}\} = \lambdatildehat_3 - \lambdatildehat_{4} = 1.1340$ as well, we can conclude that $0.2259 >  \Vert \  \Epsilonb+ \Etabtilde  \Vert_F > 0.1667$.

Overall, we covered three different outcomes. Choosing the threshold based on an OLS regression point of view validates the block containing the single variable $\Xtilde_1$. Hence, we obtain the block $\Xtilde_1$ and the block $(\Xtilde_2,\Xtilde_3,\Xtilde_4,\Xtilde_5)$ respectively. The block consisting of $\Xtilde_2$ cannot be detected by PLA with threshold $\tau_{\alpha,2}$. Therefore, $\Xtilde_2$ is not discarded by MLR, however, it is by PLA with threshold $\tau = 0.3$. Lastly, $\Xtilde_3$ is considered to be redundant in an regression-sense but not by PLA.

\section{Concluding Remarks}
\label{s:conclusion}

We have provided conditions under which PLA and MLR select the same variables and hence discard the same variables. The intuition of those conditions is that the perturbation of PLA together with the noise shall not billow more than pure noise itself. Further, we have contributed bounds for the threshold $\tau$ in PLA that ensure discarding by both PLA and regression analysis. This helps to adjust $\tau$ according to the underlying data rather than applying fixed cut-off values. We also gave bounds for the perturbation matrices.

The choice of $\tau = \tau_\alpha$ is based on the test for regression coefficients with significance level $\alpha$. With the threshold being a function of $\alpha$, $\tau_\alpha$ increases when $\alpha$ decreases. A small $\alpha$ means that the likelihood to test correctly if the regression coefficient equals zero is large. However, it increases the likelihood of conducting a type II error which is not addressed yet and should be part of future research.

Also, our derivations are based on the Pillai-Bartlett trace test which can be constructed to follow an $F$ distribution approximately. Since such approximations also exist for Roy’s largest root, Wilks’ Lambda, and Hotelling-Lawley trace, derivations based on those tests might be useful for small sample size approximations. Hotelling-Lawley trace is similarly constructed like Pillai-Bartlett trace, however, the other tests consider determinants rather than traces which makes it hard to obtain bounds with respect to matrix norms.

\appendix
\section{Complementary results for Section~\ref{s:example}}
\label{a:compl}

In the following, we provide complement material for the example in Section~\ref{s:example}. All values are rounded to four decimal places.

The sample covariance matrix $\Sigmabtildehat \equiv \Epsilonb + \Etabtilde$ of the $100 \times 5$ sample $\xbtilde = (\xbtilde_1,\ldots,\xbtilde_5)$ is given by
\begin{equation*}
\Sigmabtildehat \approx \begin{pmatrix}
8.6627 & -0.4480 & -0.6050 &  0.5438 &  0.5883 \\
-0.4480 &  1.3685 & -0.3523 & -0.9253 &  0.9582 \\
-0.6050 & -0.3523 &  5.6147 &  0.7793 & -0.2947 \\
0.5438 & -0.9253 &  0.7793 &  7.2061 & -2.0013 \\
0.5883 &  0.9582 & -0.2947 & -2.0013 &  8.4064
\end{pmatrix}
\end{equation*}
with inverse
\begin{equation*}
\Sigmabtildehat^{-1} \approx \begin{pmatrix}
  0.1207 &   0.0486 &   0.0165 &  -0.0090 &  -0.0155 \\
  0.0486 &   0.8698 &   0.0443 &   0.0805 &  -0.0818 \\
  0.0165 &   0.0443 &   0.1848 &  -0.0166 &  -0.0037 \\
-0.0090 &   0.0805 &  -0.0166 &   0.1596 &   0.0289 \\
-0.0155 &  -0.0818 &  -0.0037 &   0.0289 &   0.1361
\end{pmatrix} \; .
\end{equation*}
The estimated regression coefficients for the multivariate linear regression model \linebreak $(\xbtilde_1,\xbtilde_2,\xbtilde_4,\xbtilde_5) = \xbtilde_3 \betabtilde + \ubtilde$ are given by $\betabtildehat \approx (-0.1079, -0.0627, 0.1388, -0.0525)$, and the test result for the null of $\betabtilde = \bm{0}$ can be found in Table~\ref{table:testresults}.

\begin{table}[h!]
\centering
\caption{Test results for Roy's largest root, Wilks' Lambda, Hotelling-Lawley trace, and Pillai-Bartlett trace with respective approximated $F$-test statistics for the null of $\betabtilde = \bm{0}$.}
\label{table:testresults}
\begin{tabular}{llll}
\toprule
 & test statistic & approximated $F$ & $p$-value  \\
  \cmidrule(l{3pt}r{3pt}){1-1} \cmidrule(l{3pt}r{3pt}){2-4}
Roy 						& 0.0374 & 0.8978 & 0.4685\\
Wilks 					& 0.9639 & 0.8978 & 0.4685\\
Hotelling-Lawley	& 0.0374 & 0.8978 & 0.4685\\
Pillai-Bartlett		& 0.0361 & 0.8978 & 0.4685\\
\toprule
\end{tabular}
\end{table}

\section{Proofs}

We give the proofs for our work in this section. If terms bounded by at least $N^{-1}$ are approximated to zero, we indicate it using the notation "$\cong 0$". For example, we write $\eta \cong 0$ if $\eta = \bigO(N^{-1})$ is approximated to zero.

\begin{proof}[Proof of Lemma~\ref{th:intuition}]
Firstly, we recap approximations for the inverse sample covariance matrices (see e.g. \cite{SS90}): %p. 130 and p.118 Theorem 2.5: consistent matrix norm (hence we use Frobenius norm, since we are using the eucldean vector norm)

\begin{align}
\Sigmabhat^{-1} = \Sigmab^{-1} - \Sigmab^{-1} \Etab \Sigmab^{-1} + \bigO(\Vert\Etab\Vert_F^2) \cong \Sigmab^{-1} - \Sigmab^{-1} \Etab \Sigmab^{-1} \label{eq:approxSigmahat} 
\end{align}
and
\begin{align}
\Sigmabtildehat^{-1} & = \Sigmab^{-1} - \Sigmab^{-1} (\Epsilonb+\Etabtilde) \Sigmab^{-1} + \bigO(\Vert\Epsilonb+\Etabtilde\Vert_F^2) \nonumber \\ 
& \cong \Sigmab^{-1} - \Sigmab^{-1} (\Epsilonb+\Etabtilde)  \Sigmab^{-1}  +  \bigO(\Vert\Epsilonb\Vert_F^2)   .  \label{eq:approxSigmatildehat} 
\end{align}

The approximation error is bounded by $\bigO\left(N^{-1}\right)$ since both, $\Etab$ and $\Etabtilde$ are $\bigO(N^{-1/2})$ which is shown in Lemma~1 of \citet{NW90}. Due to (\ref{eq:approxSigmahat}), (\ref{eq:approxSigmatildehat}), and Assumption~\ref{ass:independence}, \ref{ass:meancentred} and \ref{ass:OLS}, the OLS estimators for the the linear models $\xb_\bK = \xb_\bD \betab + \ub$ and $\xbtilde_\bK = \xbtilde_\bD \betabtilde + \ubtilde$ are approximately given by
\begin{align*}
\betabhat \equiv \Sigmabhat_2^{-1} \Sigmabhat_{12}^\top \cong (\Sigmab_2^{-1} - \Sigmab_2^{-1} \Etab_2 \Sigmab_2^{-1} ) \Sigmabhat_{12}^\top \equiv (\Sigmab_2^{-1} - \Sigmab_2^{-1} \Etab_2 \Sigmab_2^{-1} ) \Etab_{12}^\top \cong \Sigmab_2^{-1} \Etab_{12}^\top
\end{align*}
and
\begin{align*}
\betabtildehat & \equiv \Sigmabtildehat_2^{-1} \Sigmabtildehat_{12}^\top \\
& \cong  \left( \Sigmab_2^{-1} - \Sigmab_2^{-1} (\Epsilonb_2+\Etabtilde_2)  \Sigmab_2^{-1}  +  \bigO(\Vert\Epsilonb\Vert_F^2) \right) \Sigmabtildehat_{12}^\top\\
& =  \left( \Sigmab_2^{-1} - \Sigmab_2^{-1} (\Epsilonb_2+\Etabtilde_2)  \Sigmab_2^{-1}  +  \bigO(\Vert\Epsilonb\Vert_F^2) \right) (\Epsilonb_{12}^\top + \Etabtilde_{12}^\top) \\
\begin{split}
& \cong \Sigmab_2^{-1}  (\Epsilonb_{12}^\top + \Etabtilde_{12}^\top)  - \Sigmab_2^{-1}\Epsilonb_2 \Sigmab_2^{-2} \Epsilonb_{12}^\top + \Sigmab_2^{-1}\Epsilonb_2 \Sigmab_2^{-2} \Etabtilde_{12}^\top + \Sigmab_2^{-1} \Etabtilde_2 \Sigmab_2^{-1} \Epsilonb_{12}^\top \\
& \qquad\qquad\qquad + \bigO(\Vert\Epsilonb\Vert_F^2) \cdot (\Epsilonb_{12}^\top + \Etabtilde_{12}^\top) \; .
\end{split}
\end{align*}

The $\{\betahat_k\}_{k\in\bK}$ are discarded according to OLS by assumption of Lemma~\ref{th:intuition}. Hence, the $\{\betatildehat_k\}_{k\in\bK}$ are discarded if $|\betatildehat_k^{(d)}| \leq |\betahat_k^{(d)}|$ for all $(k,d) \in \bK \times \bD$. However, because $\Etab$ as well as $\Etabtilde$ are $\bigO(N^{-1/2})$ we can conclude using the triangle inequality 
\begin{align*}\begin{split}
|\betatildehat_k^{(d)}| -|\betahat_k^{(d)}| & \leq | \Sigmab_2^{-1}  (\Epsilonb_{12}^\top + \Etabtilde_{12}^\top)  | + | \Sigmab_2^{-1}\Epsilonb_2 \Sigmab_2^{-2} \Epsilonb_{12}^\top | + | \Sigmab_2^{-1}\Epsilonb_2 \Sigmab_2^{-2} \Etabtilde_{12}^\top  | \\
& \qquad+ | \Sigmab_2^{-1} \Etabtilde_2 \Sigmab_2^{-1} \Epsilonb_{12}^\top | + | \bigO(\Vert\Epsilonb\Vert_F^2) \cdot (\Epsilonb_{12}^\top + \Etabtilde_{12}^\top) |  - \Sigmab_2^{-1}  \Etab_{12}^\top\\
& \leq 0
\end{split}
\end{align*}
that $\Epsilonb_{12} = \bigO(N^{-1/2})$ and $\Epsilonb_2 = \bigO(N^{-1/2})$. Hence,
\begin{align*}
\betabtildehat  & \cong \Sigmab_2^{-1}  (\Epsilonb_{12}^\top + \Etabtilde_{12}^\top)    \;. 
\end{align*}

For any $\dstar \in \bD$ with $\betahat_k^{(\dstar)} \geq 0$ the condition $|\betatildehat_k^{(\dstar)}| \leq |\betahat_k^{(\dstar)}| $ is satisfied if $ \betahat_k^{(\dstar)} \geq \betatildehat_k^{(\dstar)} \geq -\betahat_k^{(\dstar)} $. Since $\betahat_k^{(\dstar)} \cong \left(  \Sigmab_2^{-1}  \Etab_{12}^\top \right)_{[\{\dstar\},\{k\}]} = (\Sigmab_2^{-1})_{[\{\dstar\},\bD]} {\Etab_{12}^\top}_{[\bD,\{k\}]} $ and $\betatildehat_k^{(\dstar)} \cong \left(  \Sigmab_2^{-1}  (\Epsilonb_{12}^\top + \Etabtilde_{12}^\top)\right)_{[\{\dstar\},\{k\}]} = (\Sigmab_2^{-1})_{[\{\dstar\},\bD]} (\Epsilonb_{12}^\top + \Etabtilde_{12}^\top)_{[\bD,\{k\}]}   $ we obtain the desired result. The case that $\betahat^{(\dstar)} \leq 0$ remains analogously. 
\end{proof}

\begin{proof}[Proof of Theorem~\ref{l:bounds}]
Using the derivations for $\betahat_k^{(d)}$ and $\betatildehat_k^{(d)}$ from the proof of \linebreak Lemma~\ref{th:intuition}, we see that $|\betahat_k^{(d)}| \cong  | (\Sigmab_2^{-1})_{[\{\dstar\},\bD]}  (\Etab_{12}^\top)_{[\bD,\{k\}]} |\Vert (\Sigmab_2^{-1})_{[\{\dstar\},\bD]}  \Vert_2 \Vert (\Sigmab_2)^{-1}_{[\{\dstar\},\bD]}\Vert_2^{-1} $  \linebreak  $  | (\Sigmab_2^{-1})_{[\{\dstar\},\bD]}  (\Etab_{12}^\top)_{[\bD,\{k\}]} |$. Further, $|\betatildehat^{(d)}| \cong | (\Sigmab_2^{-1})_{[\{\dstar\},\bD]} ( \Epsilonb_{12}^\top+ \Etabtilde_{12}^\top)_{[\bD,\{k\}]} | \leq$     \linebreak       $ \Vert (\Sigmab_2^{-1})_{[\{\dstar\},\bD]}   \Vert_2 \Vert ( \Epsilonb_{12}^\top+ \Etabtilde_{12}^\top)_{[\bD,\{k\}]}  \Vert_2$ by Hölder's inequality. Hence, the condition that $|\betahat_k^{(d)}| \leq |\betatildehat_k^{(d)}|$ is satisfied if  $\Vert (\Sigmab_2^{-1})_{[\{\dstar\},\bD]}  \Vert_2  \Vert (\Sigmab_2^{-1})_{[\{\dstar\},\bD]}\Vert_2^{-1}   | (\Sigmab_2^{-1})_{[\{\dstar\},\bD]}  (\Etab_{12}^\top)_{[\bD,\{k\}]} | \leq \Vert (\Sigmab_2^{-1})_{[\{\dstar\},\bD]}   \Vert_2 \Vert ( \Epsilonb_{12}^\top+ \Etabtilde_{12}^\top)_{[\bD,\{k\}]}  \Vert_2$. The result follows by shortening $\Vert (\Sigmab_2^{-1})_{[\{\dstar\},\bD]}    \Vert_2 $.
\end{proof}
\begin{proof}[Proof of Corollary~\ref{c:bounds}]
The corollary is an immediate result of Theorem~\ref{l:bounds} since $
(\Vert (\Sigmab_2)^{-1}_{[\{\dstar\},\bD]}\Vert_2)^{-1}   |  (\Sigmab_2^{-1})_{[\{\dstar\},\bD]}  (\Etab_{12}^\top)_{[\bD,\{k\}]} | \leq \Vert (\Etab_{12}^\top)_{[\bD,\{k\}]}  \Vert_2$ by Hölder's inequality.
\end{proof}
\begin{proof}[Proof of Corollary~\ref{c:convergencerate}]
That $\Epsilonb_2$ and $\Epsilonb_{12}$ are $\bigO\left( N^{-1/2}\right)$ follows from the proof of Lemma~\ref{th:intuition}. Reproducing the proof of Lemma~\ref{th:intuition} using the MLR models $\xb_\bD = \xb_\bK \bm{\gamma} + \ub$ and $\xbtilde_\bD = \xbtilde_\bK \tilde{\bm{\gamma}} + \ubtilde$ respectively yields that $\Epsilonb_1 = \bigO\left( N^{-1/2}\right)$ as well.
\end{proof}
\begin{proof}[Proof of Theorem~\ref{l:tau}]
Our derivations are based on the Pillai-Bartlett trace test which follows under the null of $\betab = \bm{0}$ the test statistic 
\begin{equation*}
{\rm PB} \equiv  \tr (  (\Sigmabtildehat - \Sigmabtildehat_{\ub} ) \Sigmabtildehat_{\ub}^{-1} ) 
\end{equation*}
where $(N-1) \Sigmabtildehat_{\ub} \equiv (\xbtilde_\bK - \xbtilde_\bD \betabtildehat)^\top (\xbtilde_\bK - \xbtilde_\bD \betabtildehat) $ and $(N-1) \Sigmabtildehat \equiv \xbtilde_\bK^\top \xbtilde_\bK $ such that $(N-1) \Sigmabtildehat_{\ub}$ is the residual sum of squares and cross products matrix of the full model, and $ (N-1)  (\Sigmabtildehat - \Sigmabtildehat_{\ub} ) $ is the sum of squares and cross products matrix under the null that is to say under the reduced model. Under the null of $\betab = \bm{0}$, it holds that 
\begin{equation*}
F_{{\rm PB} }\dfrac{{\rm PB} \cdot D_2 }{(S-{\rm PB}) D_1 } \sim F_{D_1,D_2}
\end{equation*}
where $D_1 \equiv D  K$, $D_2 \equiv S(N+S-K)$. and $S\equiv \min\{D,K\}$. Therefore, the $\{\Xtilde_d\}$ are discarded in a MLR sense if $F_{{\rm PB} } < F_\alpha$ where $F_\alpha$ is the quantile of the $F$-distribution with degrees of freedom $D_1 \equiv D  K$ and $D_2 \equiv S(N+S-K)$ (\citet{JW02}, \citet{MKB79} and \citet{MP84}).

By definition, $(N-1) \Sigmabtildehat_{\ub} \equiv \xbtilde_\bK^\top \xbtilde_\bK \equiv (N-1) \Sigmabtildehat_1$. Further, it holds that $(N-1)(\Sigmabtildehat - \Sigmabtildehat_{\ub}) \equiv \xbtilde_\bK^\top \xbtilde_\bK - (\xbtilde_\bK - \xbtilde_\bD \betabtildehat)^\top (\xbtilde_\bK - \xbtilde_\bD \betabtildehat) \equiv (N-1)  (\Epsilonb_{12} + \Etabtilde_{12}) \Sigmabtildehat_2^{-1} (\Epsilonb_{12} + \Etabtilde_{12})^\top  $ which can be obtained by matrix algebra using $\betabtildehat \equiv (\xb_\bD^\top \xb_\bD)^{-1} \xb_\bD^\top \xb_\bK$.

$\Sigmabtildehat_{\ub}$ is positive definite by construction and $\Sigmabtildehat - \Sigmabtildehat_{\ub}$ is positive semi-definite because $\Sigmabtildehat_2^{-1}$ is positive definite. Hence
\begin{align*}
\tr ( (\Sigmabtildehat - \Sigmabtildehat_{\ub}) \Sigmabtildehat_{\ub}^{-1} )^2 & \leq \tr ( (\Sigmabtildehat - \Sigmabtildehat_{\ub})^2 ) \tr ( (\Sigmabtildehat_{\ub}^{-1})^2 ) \\
& = \Vert \Sigmabtildehat - \Sigmabtildehat_{\ub} \Vert_F \Vert \Sigmabtildehat_{\ub}^{-1} \Vert_F \\
 & \equiv \Vert (\Epsilonb_{12} + \Etabtilde_{12}) \Sigmabtildehat_2^{-1} (\Epsilonb_{12} + \Etabtilde_{12})^\top \Vert_F \Vert  \Sigmabtildehat_1^{-1} \Vert_F \\
 & \leq \Vert \Epsilonb + \Etabtilde \Vert_F^2 \Vert \Sigmabtildehat_1^{-1} \Vert_F \Vert \Sigmabtildehat_2^{-1} \Vert_F 
\end{align*}
and we can conclude that $\tr ( (\Sigmabtildehat - \Sigmabtildehat_{\ub}) \Sigmabtildehat_{\ub}^{-1}  ) = [ \tr( (\Sigmabtildehat - \Sigmabtildehat_{\ub}) \Sigmabtildehat_{\ub}^{-1} )^2 ]^{1/2} \leq$ \linebreak $  \Vert \Epsilonb + \Etabtilde \Vert_F \Vert \Sigmabtildehat_1^{-1} \Vert_F^{1/2} \Vert \Sigmabtildehat_2^{-1} \Vert_F^{1/2}$.

Due to Theorem 2 of \citet{BD21} which is based on Corollary 1 of \citet{YW15}, it holds that the variables $\{\Xtilde_d\}$ are discarded both by PLA and OLS if for all $\delta \in  \bDelta$
\begin{equation*}
\dfrac{2^{2/3}\Vert \  \Epsilonb+ \Etabtilde  \Vert_F }
{\min\{\lambda_{\delta-1} - \lambda_\delta, \lambda_\delta - \lambda_{\delta+1}\}} 
\leq \tau \leq 
 \dfrac{2^{2/3}  F_\alpha \Vert \Sigmabtildehat_1^{-1} \Vert_F^{-1/2} \Vert \Sigmabtildehat_2^{-1} \Vert_F^{-1/2}  S  D_1  D_2^{-1} }
 {\min\{\lambda_{\delta-1} - \lambda_\delta, \lambda_\delta - \lambda_{\delta+1}\}}  \;.
\end{equation*}
All $D$ equations hold if the condition is fulfilled for the smallest eigengap $ \min_{\delta\in\bDelta} \{\lambda_{\delta-1} - \lambda_\delta, \lambda_\delta - \lambda_{\delta+1}\} $. 
\end{proof}
\begin{proof}[Proof of Corollary~\ref{c:perturbationbounds}]
The result follows from the proof of Theorem~\ref{l:tau} because we fail to reject the null that $\betabtilde = \bm{0}$ if $\Vert \  \Epsilonb+ \Etabtilde  \Vert_F > 2^{2/3}   \Vert \Sigmabtildehat_1^{-1} \Vert_F^{-1/2} \Vert \Sigmabtildehat_2^{-1} \Vert_F^{-1/2}  \dfrac{S D_1}{F_\alpha^{-1}D_2 + D_1}$
\end{proof}
\begin{proof}[Proof of Corollary~\ref{c:tau}]
The corollary is an immediate result of Theorem~\ref{l:tau} and Corollary~\ref{c:bounds} since $
\Vert (\Epsilonb_{12}^\top + \Etabtilde_{12}^\top)_{[\bD,\{k\}]}  \Vert_2 \leq \Vert (\Epsilonb + \Etabtilde) \Vert_F $ and since $\Vert (\Etab_{12}^\top)_{[\bD,\{k\}]} \Vert_2 \leq \Vert (\Etab_{12}^\top)_{[\bD,\bK]} \Vert_F$ for all $k\in\bK$.
\end{proof}
\begin{proof}[Proof of Corollary~\ref{c:tau2}]
The corollary is an immediate result of Theorem~\ref{l:tau}  since $ \min_{\delta =\deltastar} \{ \lambda_{\delta-1} - \lambda_\delta, \lambda_\delta - \lambda_{\delta+1} \} = \min\{\lambda_{\deltastar-1} - \lambda_\deltastar, \lambda_\deltastar - \lambda_{\deltastar+1} \}$.
\end{proof}
\begin{proof}[Proof of Corollary~\ref{c:tau3}]
The corollary is an immediate result of Corollary~\ref{c:tau} since $ \min_{\delta =\deltastar} \{ \lambda_{\delta-1} - \lambda_\delta, \lambda_\delta - \lambda_{\delta+1} \} = \min\{\lambda_{\deltastar-1} - \lambda_\deltastar, \lambda_\deltastar - \lambda_{\deltastar+1} \}$.
\end{proof}

\section*{References}
\bibliographystyle{elsarticle-harv} 
\bibliography{bib}{}

\end{document}